\documentclass[graybox]{svmult}
\usepackage{mathptmx}
\usepackage{helvet}
\usepackage{courier}
\usepackage{makeidx}
\usepackage{graphicx}
\usepackage{multicol}
\usepackage[bottom]{footmisc}

\usepackage{pgfplots}
\usepackage{pgfplotstable}

\usetikzlibrary{external}
\makeindex

\begin{document}
\title*{ArbiLoMod: Local Solution Spaces by Random Training in Electrodynamics}
\author{Andreas Buhr, Christian Engwer, Mario Ohlberger, and Stephan Rave}
\titlerunning{ArbiLoMod-Training for Electrodynamics}
\authorrunning{Andreas Buhr et.~al.}
\institute{
Andreas Buhr \email{andreas@andreasbuhr.de}
\and Christian Engwer \email{christian.engwer@uni-muenster.de}
\and Mario Ohlberger \email{mario.ohlberger@uni-muenster.de}
\and Stephan Rave \email{stephan.rave@uni-muenster.de}
\at Institute for Computational and Applied Mathematics, University of M\"unster, Einsteinstra\ss e 62, D-48149 M\"unster, Germany 
}
\maketitle

\abstract{
The simulation method ArbiLoMod \cite{BEOR15} has the goal of providing users of Finite Element based simulation software with quick re-simulation after localized changes to the model under consideration. It generates a Reduced Order Model (ROM) for the full model without ever solving the full model. To this end, a localized variant of the Reduced Basis method is employed, solving only small localized problems in the generation of the reduced basis. The key to quick re-simulation lies in recycling most of the localized basis vectors after a localized model change. In this publication, ArbiLoMod's local training algorithm is analyzed numerically for the non-coercive problem of time harmonic Maxwell's equations in 2D, formulated in $H(\mathrm{curl})$.
}
\section{Introduction}
\begin{figure}
\centering
\includegraphics[width=.9\textwidth,viewport=0 220 1476 470, clip=true]{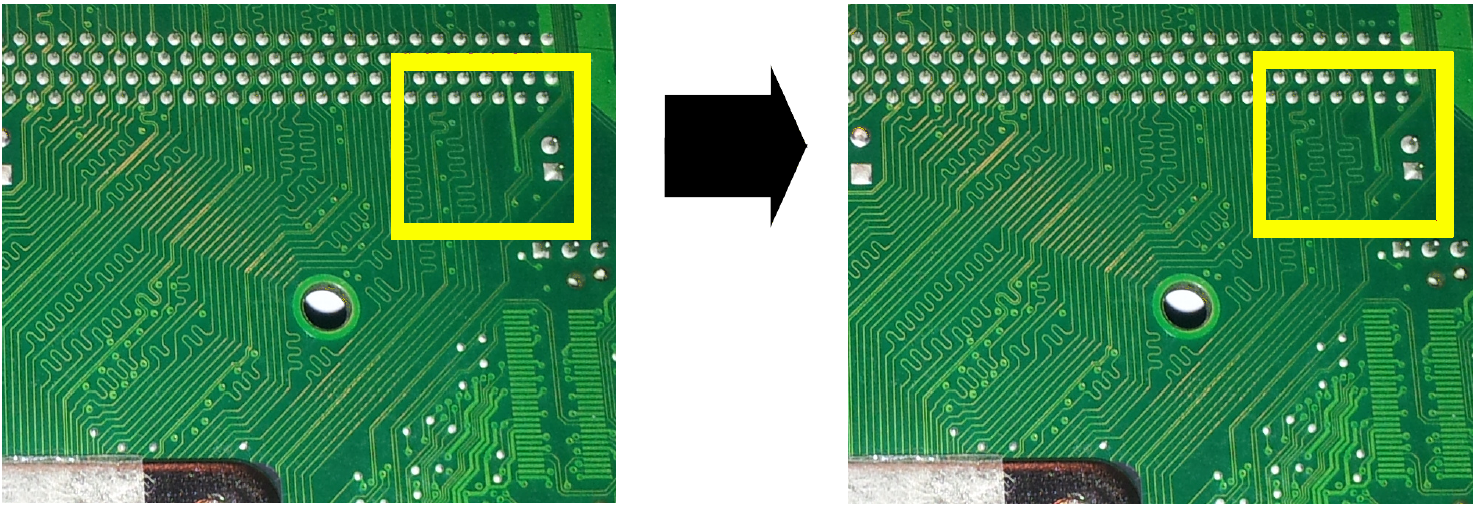}
\caption{Printed circuit board subject to local modification of conductive tracks.}
\label{fig:boards}
\end{figure}

Simulation software based on the Finite Element Method is an essential ingredient of many engineering workflows. 
In their pursue of design goals, engineers often simulate structures several times, applying small changes
after each simulation. This results in large similarities between subsequent simulation runs. These
similarities are usually not considered by simulation software. The simulation method ArbiLoMod was
designed to change that and accelerate the subsequent simulation of geometries which only differ in small details.
A motivating example is the design of mainboards for PCs. Improvements in the signal integrity properties of
e.g. DDR memory channels is often obtained by localized changes, as depicted in Fig.~\ref{fig:boards}.

ArbiLoMod was also designed with the available computing power in mind: Today, cloud environments are just a few
clicks away and everyone can access hundreds of cores easily. However, the network connection to these cloud
environments is relatively slow in comparison to the available computing power, so a method which should perform well
under these circumstances must be designed to be communication avoiding.

At the core of ArbiLoMod lies a localized variant of the Reduced Basis Method. The Reduced Basis Method is 
a well established approach to create reduced order models (ROMs) and its application to the 
Maxwell's equations has been extensively investigated by many groups (see e.g. \cite{Chen2010, BennerHess, EdlingerSommer2014, Pomplun2010, Fares20115532}).
On the other hand, there are lots of methods with localized basis generation (e.g. 
\cite{Maday2002a,Iapichino2012,PhuongHuynh2012,Maier2014,Efendiev2013,strouboulis2000design,OS15}).
However to the authors' knowledge, only little was published on the combination of both.
In \cite{Chen2011}, 
the Reduced Basis Element Method is applied to time harmonic Maxwell's equations.

This publication evaluates ArbiLoMod's training numerically for the time harmonic Maxwell's equation.
The remainder of this article is structured as follows:
In the following Section \ref{buhr_andreas:sec:problem_setting}, the problem setting is given. 
Section \ref{buhr_andreas:sec:arbilomod_for_maxwell} outlines ArbiLoMod and highlights the specialties when
considering inf-sup stable problems in $H(\mathrm{curl})$.
Afterwards, we demonstrate ArbiLoMod's performance on a numerical example in Section \ref{buhr_andreas:sec:numerical_example}.
Finally, we conclude in Section \ref{buhr_andreas:sec:conclusion}.
\section{Problem Setting}
\label{buhr_andreas:sec:problem_setting}

We consider Maxwell's equations \cite{maxwell1861xxv} on the polygonal domain $\varOmega$. 
The material is assumed to be linear and isotropic, i.e. the electric permittivity $\varepsilon$ and the magnetic permeability $\mu$
are scalars.
On the boundary $\partial \varOmega = \varGamma_R \cup \varGamma_D$ we impose Dirichlet 
($E \times n = 0 \quad \mathrm{on} \ \varGamma_D$)
 and Robin
($H \times n = \kappa(E \times n) \times n \quad \mathrm{on} \ \varGamma_R$, \cite[eq. (1.18)]{monk2003finite})
boundary conditions with the surface impedance parameter $\kappa$.
$n$ denotes the unit outer normal of $\varOmega$. The excitation is given by a current density $\hat j$.

For the time harmonic case, this results in the following weak formulation:\\
find $u \in V:= H(\mathrm{curl})$ so that
\begin{eqnarray}
\label{buhr_andreas:eq:problem}
a(u,v;\omega) &=& f(v;\omega) \qquad \forall v \in V, \\
a(u,v;\omega) &:=& \int_\varOmega \frac 1 \mu (\nabla \times u) \cdot (\nabla \times \overline v)
                   - \varepsilon \omega ^2 (u \cdot \overline v) \mathrm d v 
                   + i \omega \kappa \int_{\varGamma_R} (u \times n) \cdot (\overline v \times n) \mathrm d S,
\nonumber\\
f(v;\omega) &:=& - i \omega \int_\varOmega (\hat j \cdot \overline v ) \nonumber
\end{eqnarray}
where $\omega$ is the angular frequency. We see $\omega$ as a parametrization to this problem. We solve in a parameter domain sampled by a finite training set $\varXi$.

We use the inner product and energy norm given by:
\begin{eqnarray}
(v,u)_V &:=& \int_\varOmega \frac 1 \mu (\nabla \times u) \cdot (\nabla \times \overline v)
                   + \varepsilon \omega_\mathrm{max} ^2 (u \cdot \overline v) \mathrm d v
                   +  \omega_\mathrm{max} \kappa \int_{\varGamma_R} (u \times n) \cdot (\overline v \times n) \mathrm d S, \nonumber
\\ \|u\|_V &:=& \sqrt{(u,u)_V}.
\end{eqnarray}
\subsection{Discretization}
We assume there is a non overlapping domain decomposition with subdomains $\varOmega_i$, $\varOmega_i \cap \varOmega_j = \emptyset$ for $i \ne j$.
For simplicity, we assume it to be rectangular.
The domain decomposition should cover the problem domain $\varOmega \subseteq \bigcup_i \varOmega_i$ but the subdomains need not resolve the domain.
This is important as we want to increase or decrease the size of $\varOmega$ between simulation runs without changing the domain decomposition.
For example, in a printed circuit board (PCB), the metal traces are often simulated as being outside of the domain. Thus, a change of the traces leads to a change of the
calculation domain.

Further we assume there is a triangulation of $\varOmega$ which resolves the domain decomposition.
We denote by $V_h$ the discrete space spanned by lowest order Nedelec ansatz functions \cite{nedelec1980} on this triangulation.
\section{ArbiLoMod for Maxwell's Equations}
\label{buhr_andreas:sec:arbilomod_for_maxwell}
The main ingredients of ArbiLoMod are (1) a localizing space decomposition, (2) localized trainings for reduced local subspace generation, (3) a localized a-posteriori error estimator and (4) a localized enrichment for basis improvement. 
In this publication, we focus on the first two steps, which are described in the following.
\subsection{Space Decomposition}
\label{buhr_andreas:sec:space_decomposition}
Localization is performed in ArbiLoMod using a direct decomposition of the ansatz space into localized subspaces. In the 2D case with Nedelec ansatz functions, there are only volume spaces $V_{\{i\}}$ associated with the subdomains $\varOmega_i$, and interface spaces $V_{\{i,j\}}$ associated with the interfaces between $\varOmega_i$ and $\varOmega_j$. The interface spaces are only associated with an interface, they are subspaces of the global function space and have support on two domains. They are not trace spaces. In higher space dimensions and/or with different ansatz functions, there may be also spaces associated with edges and nodes of the domain decomposition \cite{BEOR15}.
\begin{equation}
\label{buhr_andreas:eq:direct_decomposition}
V_h = \left( \bigoplus_i V_{\{i\}} \right) \bigoplus \left( \bigoplus_{i,j} V_{\{i,j\}} \right)
\end{equation}
The spaces $V_{\{i\}}$ are simply defined as the span of all ansatz functions having support only on $\varOmega_i$. With $\mathcal{B}$ denoting the set of all FE basis functions, we define:
\begin{equation}
V_{\{i\}} := \mathrm{span}\left\{\psi \in \mathcal{B} \ \big| \ \mathrm{supp}(\psi) \subseteq \overline \varOmega_i \right\}.
\end{equation}
The interface spaces $V_{\{i,j\}}$ are not simply the span of FE ansatz functions.
Instead, they are defined as the span of all ansatz functions on the interface plus their extension to the adjacent subdomains. The extension is done by solving for a fixed frequency $\omega'$ with Dirichlet zero boundary conditions.
The formal definition of the interface spaces is in two steps: First, we define $U_{\{i,j\}}$ as the space spanned by all ansatz function having support on both 
$\varOmega_i$ and $\varOmega_j$:
\begin{equation}
U_{\{i,j\}} := \mathrm{span} \left\{ \psi \in \mathcal{B} \ \big| \ \mathrm{supp}(\psi) \cap \varOmega_i \ne \emptyset, \mathrm{supp}(\psi) \cap \varOmega_j \ne \emptyset \right\}.
\end{equation}
Then we define the extension operator:
\begin{eqnarray}
\mathrm{Extend} : U_{\{i,j\}} &\rightarrow& V_{\{i\}} \oplus U_{\{i,j\}} \oplus V_{\{j\}},\\
\varphi  &\mapsto& \varphi + \psi \nonumber\\
&&\mbox{where } \psi \in V_{\{i\}} \oplus V_{\{j\}} \mbox{ solves} \nonumber\\
&&
 a(\varphi + \psi, \phi; \omega') = 0 \qquad \forall \phi \in V_{\{i\}} \oplus V_{\{j\}}. \nonumber
\end{eqnarray}
We then can define the interface spaces as
\begin{equation}
V_{\{i,j\}} := \left\{ \mathrm{Extend}(\varphi) \ \big| \ \varphi \in U_{\{i,j\}} \right\}.
\end{equation}
Equation (\ref{buhr_andreas:eq:direct_decomposition}) holds for this decomposition, i.e.~there is a unique decomposition of every element of $V_h$ into the localized
subspaces. We define projection operators $P_{\{i\}} : V_h \rightarrow V_{\{i\}}$
and $P_{\{i,j\}} : V_h \rightarrow V_{\{i,j\}}$ by the relation
\begin{equation}
\varphi = \sum_i P_{\{i\}}(\varphi) + \sum_{i,j} P_{\{i,j\}}(\varphi) \qquad \forall \varphi \in V_h.
\end{equation}
\subsection{Training}
\label{buhr_andreas:sec:training}
Having defined the localized spaces, we create reduced localized subspaces
$\widetilde V_{\{i\}} \subseteq V_{\{i\}}$ and
$\widetilde V_{\{i,j\}} \subseteq V_{\{i,j\}}$ by a localized training procedure.
The training is inspired by the ``Empirical Port Reduction'' introduced by 
Eftang et.al.\cite{Eftang2013}.
Its main four steps are:
\begin{enumerate}
\item solve the problem (\ref{buhr_andreas:eq:problem}) on a small training domain around the space in question
with zero boundary values for all parameters in the training set $\varXi$,
\item
solve the homogeneous equation
repeatedly on a small training domain around the space in question with
random boundary values for all parameters in $\varXi$,
\item
apply the space decomposition to all computed local solutions
to obtain the part belonging to the space in question
and
\item
use a greedy procedure to create a space approximating this set.
\end{enumerate}
For further details, we refer to \cite{BEOR15}. The small training domain for an interface space consists of the six subdomains around that interface. The small training domain for a volume space consists of nine subdomains: the subdomain in question and the eight subdomains surrounding it.

While the ``Empirical Port Reduction'' in \cite{Eftang2013}
generates an interface space and requires ports which do not intersect,
this training can be used for both interface and volume spaces. It can
handle touching ports and can thus be applied to a standard domain decomposition.
\subsection{Reduced Model}
In these first experiments the reduced global problem is obtained by a simple Galerkin projection onto
the direct sum of all reduced local subspaces. The global solution space is
\begin{equation}
\widetilde V_h := 
\left( \bigoplus_i \widetilde V_{\{i\}} \right) \bigoplus \left( \bigoplus_{i,j} \widetilde V_{\{i,j\}} \right).
\end{equation}
And the reduced problem reads: find $\widetilde u \in \widetilde V_h$ such that
\begin{equation}
a(\widetilde u, v; \mu) = f(v) \qquad \forall v \in \widetilde V_h\ .
\end{equation}
\section{Numerical Example}
\label{buhr_andreas:sec:numerical_example}
The numerical experiments are performed with pyMOR \cite{pymor}. The source code used to reproduce the results in this publication is provided
alongside with this publication and can be downloaded at http://www.arbilomod.org/morepas2015.tgz. See the README file therein for installation instructions.
\subsection{Geometry Simulated}
\begin{figure}
\begin{center}
\includegraphics[width=0.3\textwidth]{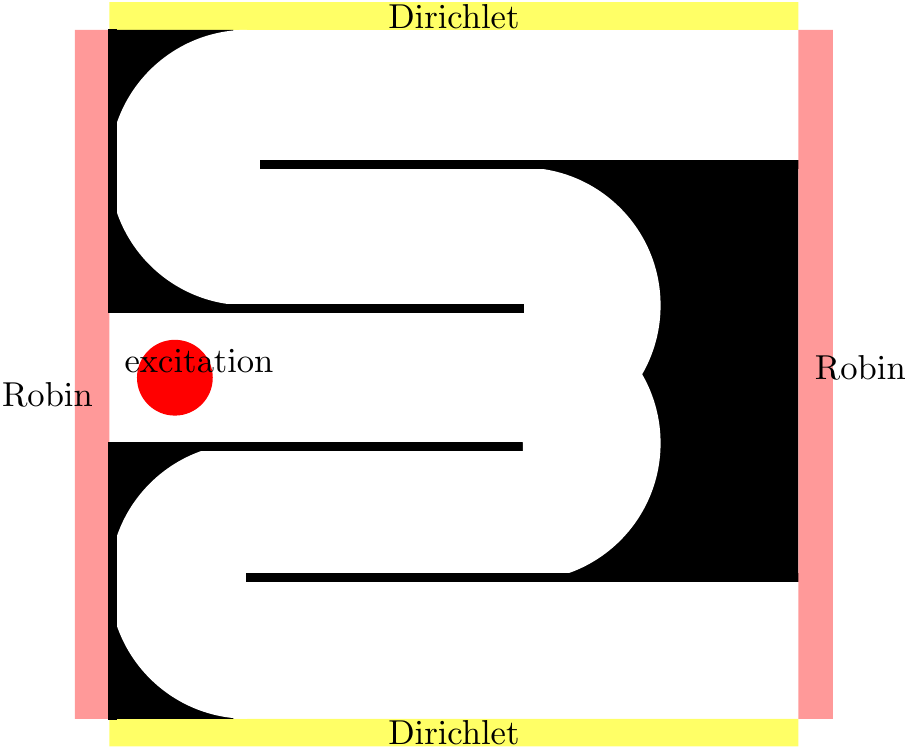}
\hspace{20pt}
\includegraphics[width=0.3\textwidth]{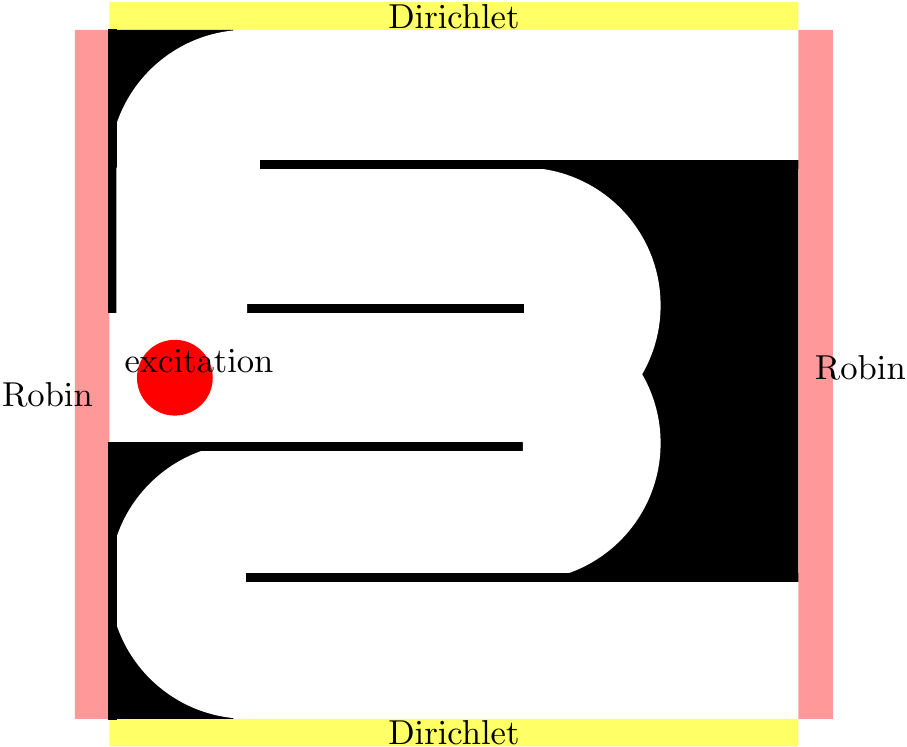}
\end{center}
\caption{Geometries simulated. Black area is not part of the domain and treated as Dirichlet zero boundary. Note the change is topology changing.}
\label{buhr_andreas:fig:example}
\end{figure}
We simulate the unit square $(0,1)\times (0,1)$ with robin boundary conditions at $\varGamma_R := 0 \times (0,1) \cup 1 \times (0,1)$ and Dirichlet zero boundary conditions at $\varGamma_D :=(0,1) \times 0 \cup (0,1) \times 1$. The surface impedance parameter $\kappa$ is chosen as the impedance of free space, $\kappa = 1 / 376.73$ Ohm. We introduce some structure by inserting perfect electric conductors (PEC) into the domain, see Fig.~\ref{buhr_andreas:fig:example}. The PEC is modeled as Dirichlet zero boundary condition. 
Note that it is slightly asymmetric intentionally, to produce more interesting behavior. The mesh does not resolve the geometry. Rather, we use a structured mesh and remove all degrees of freedom which are associated with an edge whose center is inside of the PEC structure. The structured mesh consists of 100 times 100 squares, each of which is divided into four triangles. With each edge, one degree of freedom is associated, which results in 60200 degrees of freedom, some of which are ``disabled'' because they are in PEC or on a Dirichlet boundary. The parameter domain is the range from 10 MHz to 1 GHz. For the training set $\varXi$, we use 100 equidistant points in this range, including the endpoints. To simulate an ``arbitrary local modification'', the part of the PEC within $(0.01, 0.2) \times (0.58,0.80)$ is removed and the simulation domain is enlarged.

The excitation is a current
\begin{equation}
j(x,y) := \mathrm{exp}\left( - \frac{(x - 0.1)^2 + (y - 0.5)^2}{1.25 \cdot 10^{-3}}\right) \cdot \mathrm{e}_y
\end{equation}
\begin{figure}
\begin{center}
\begin{tabular}{cccccc}
\includegraphics[width=0.15\textwidth]{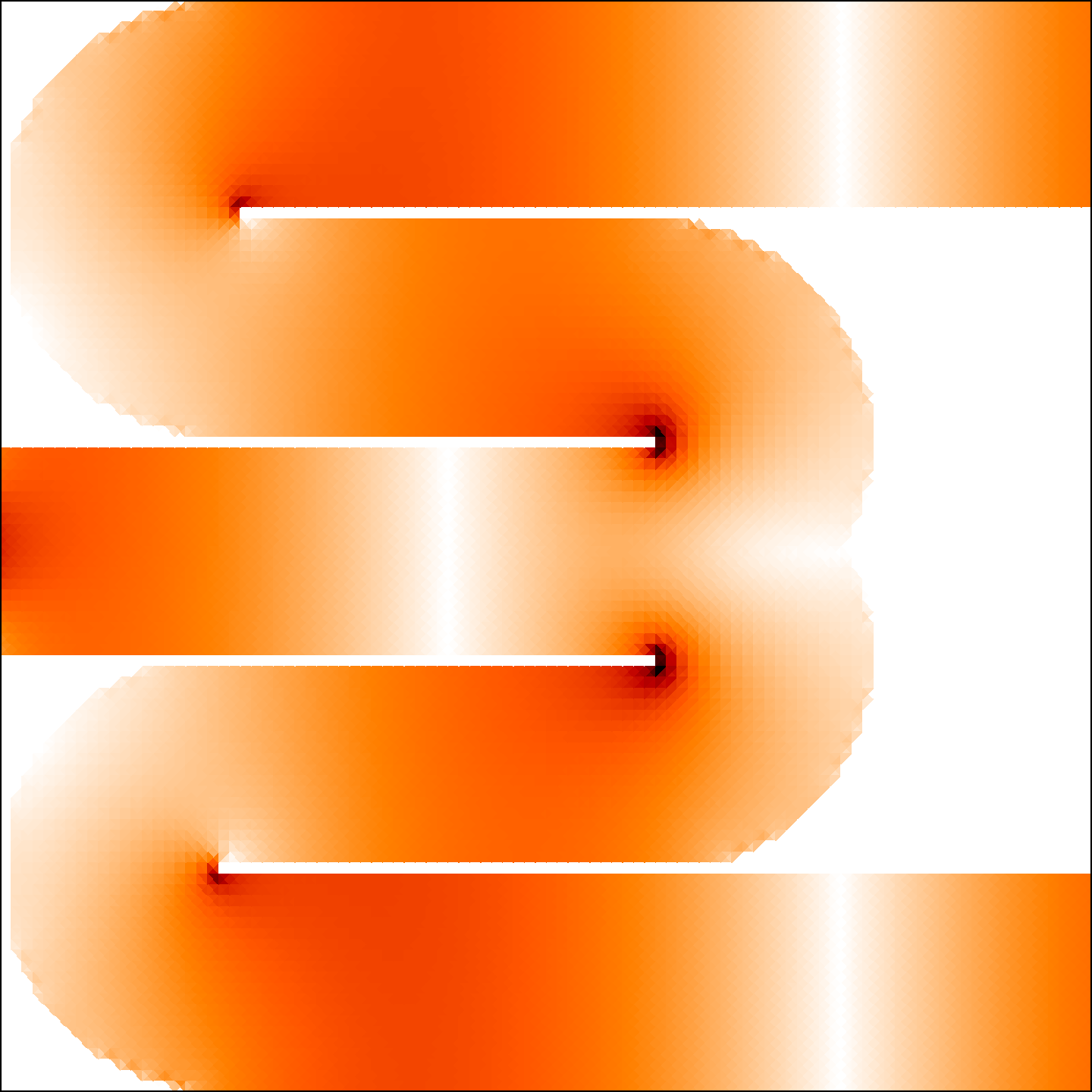}&
\includegraphics[width=0.15\textwidth]{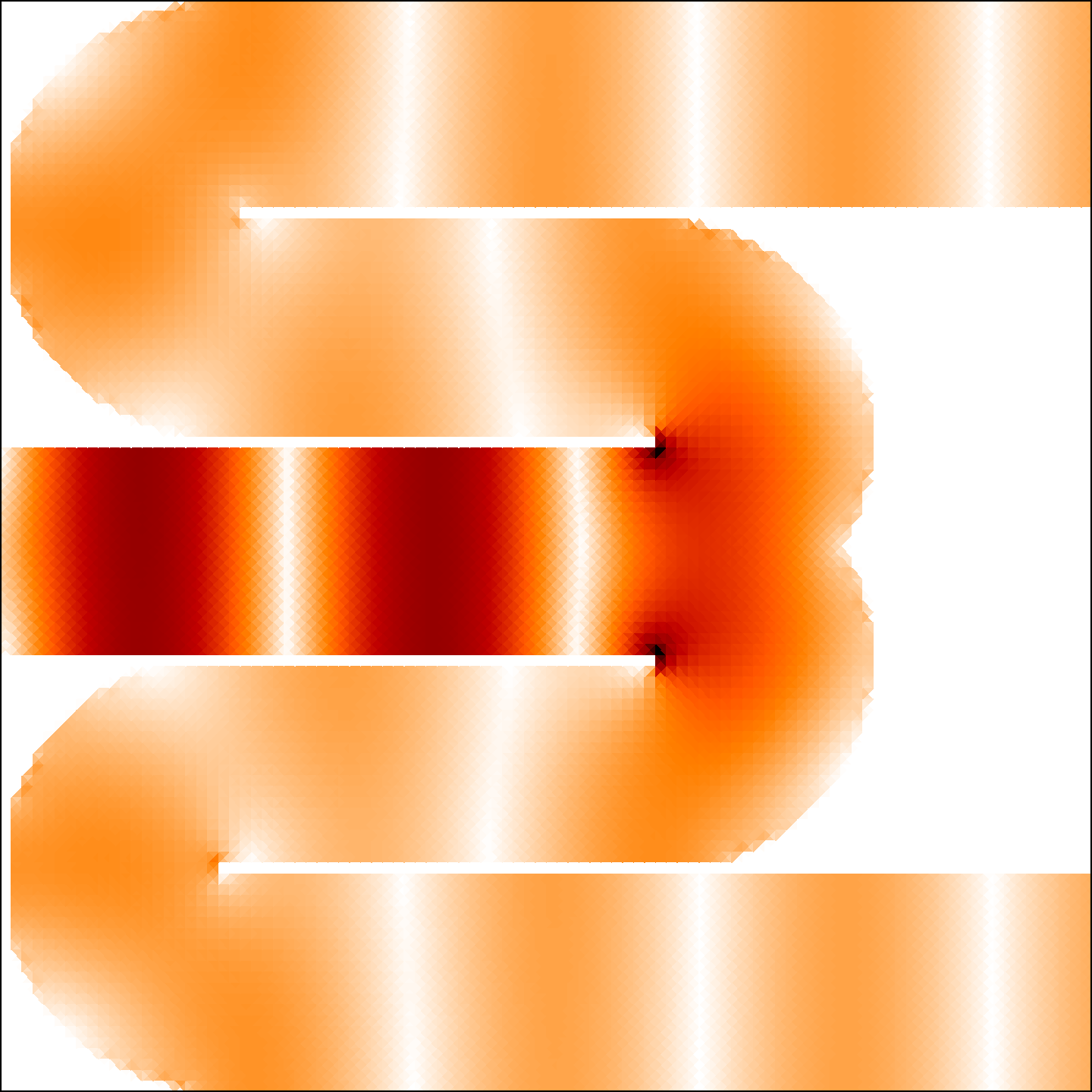}&
\includegraphics[width=0.15\textwidth]{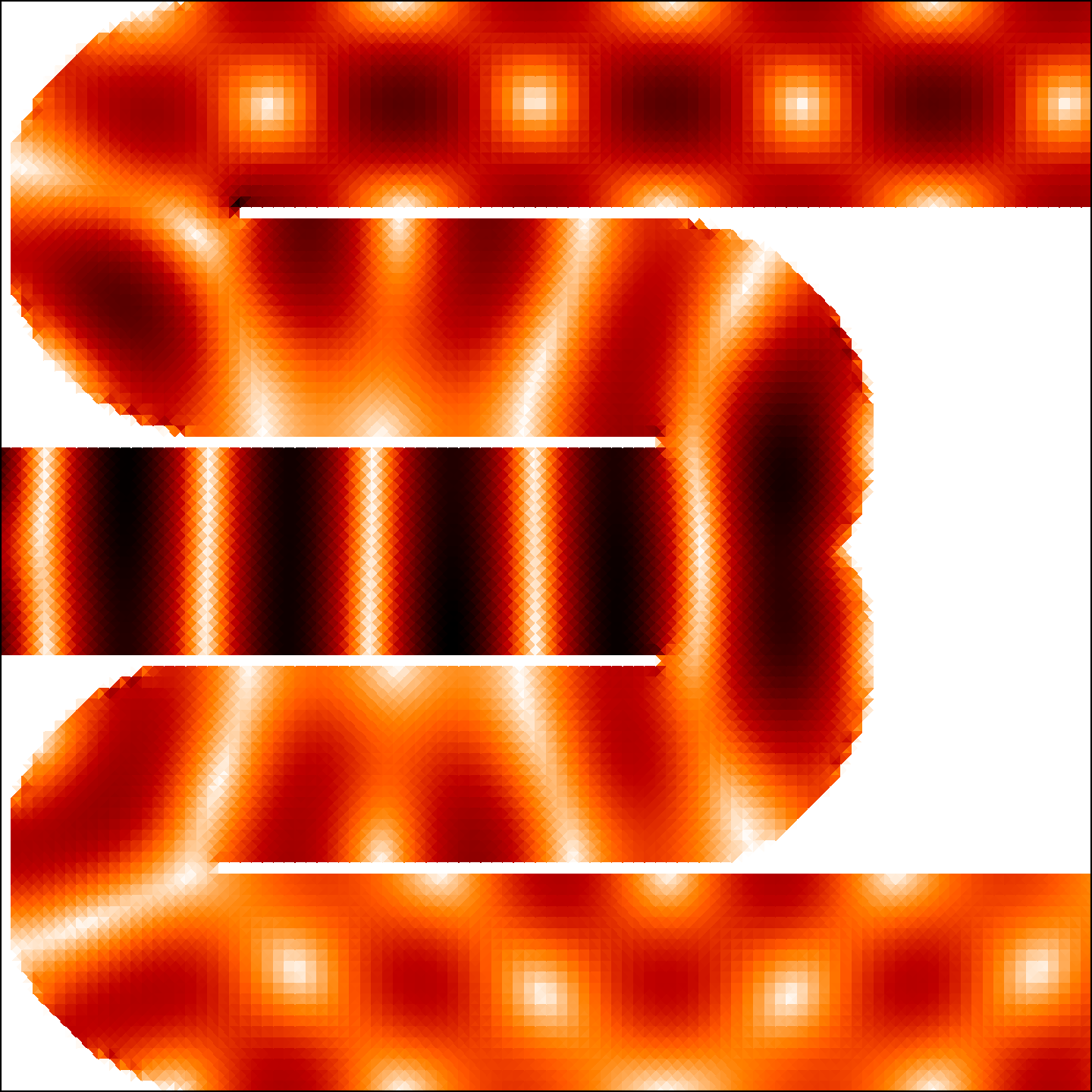}&
\includegraphics[width=0.15\textwidth]{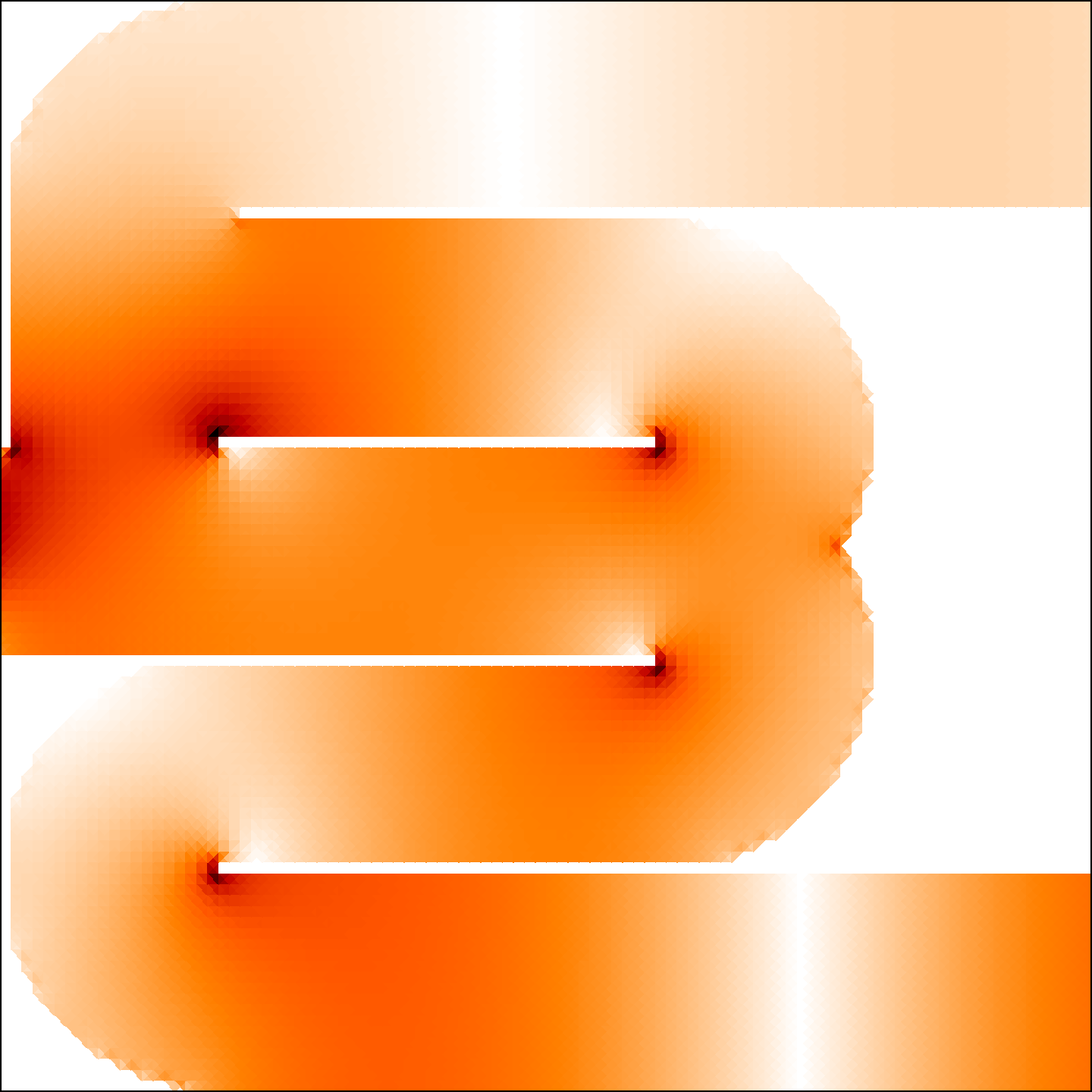}&
\includegraphics[width=0.15\textwidth]{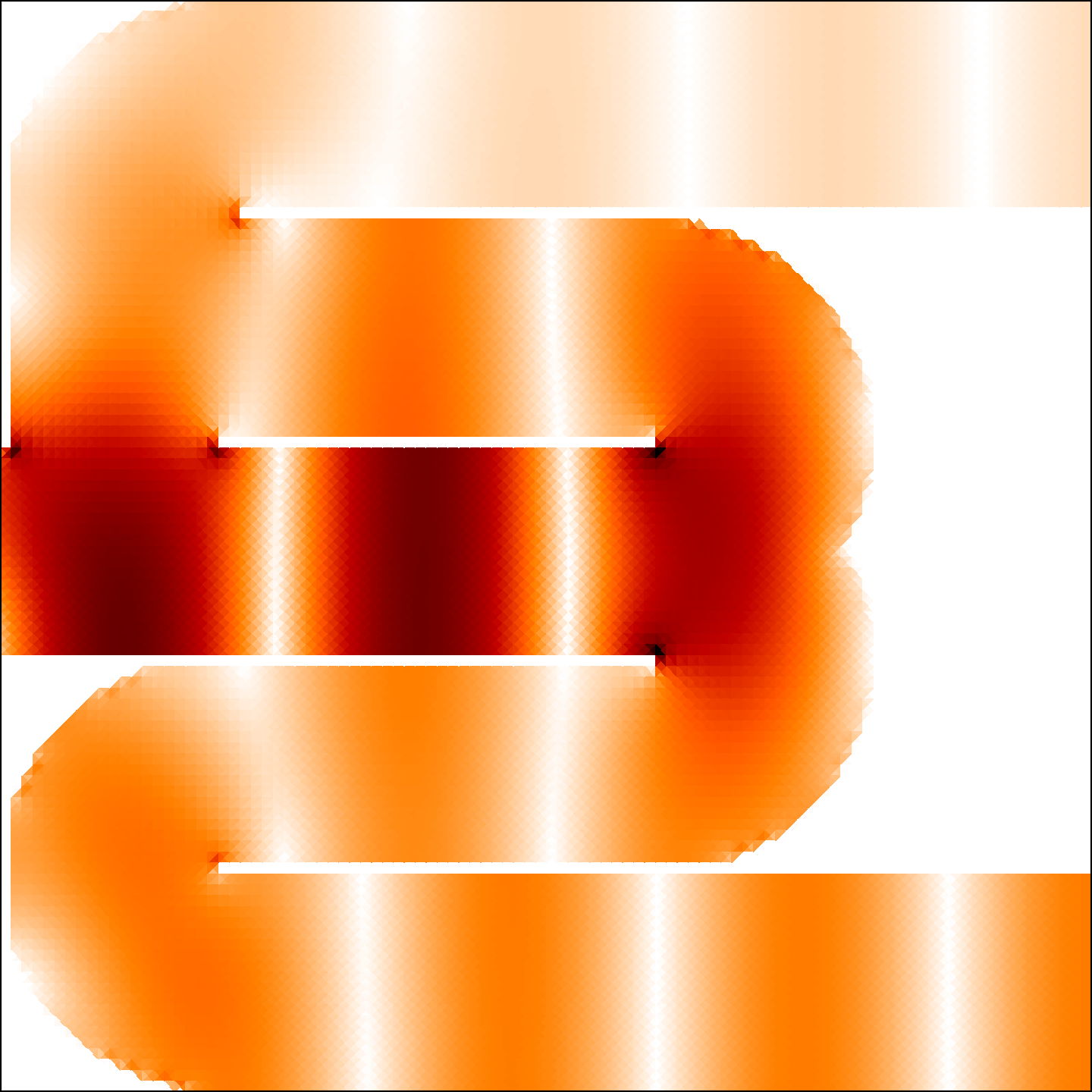}&
\includegraphics[width=0.15\textwidth]{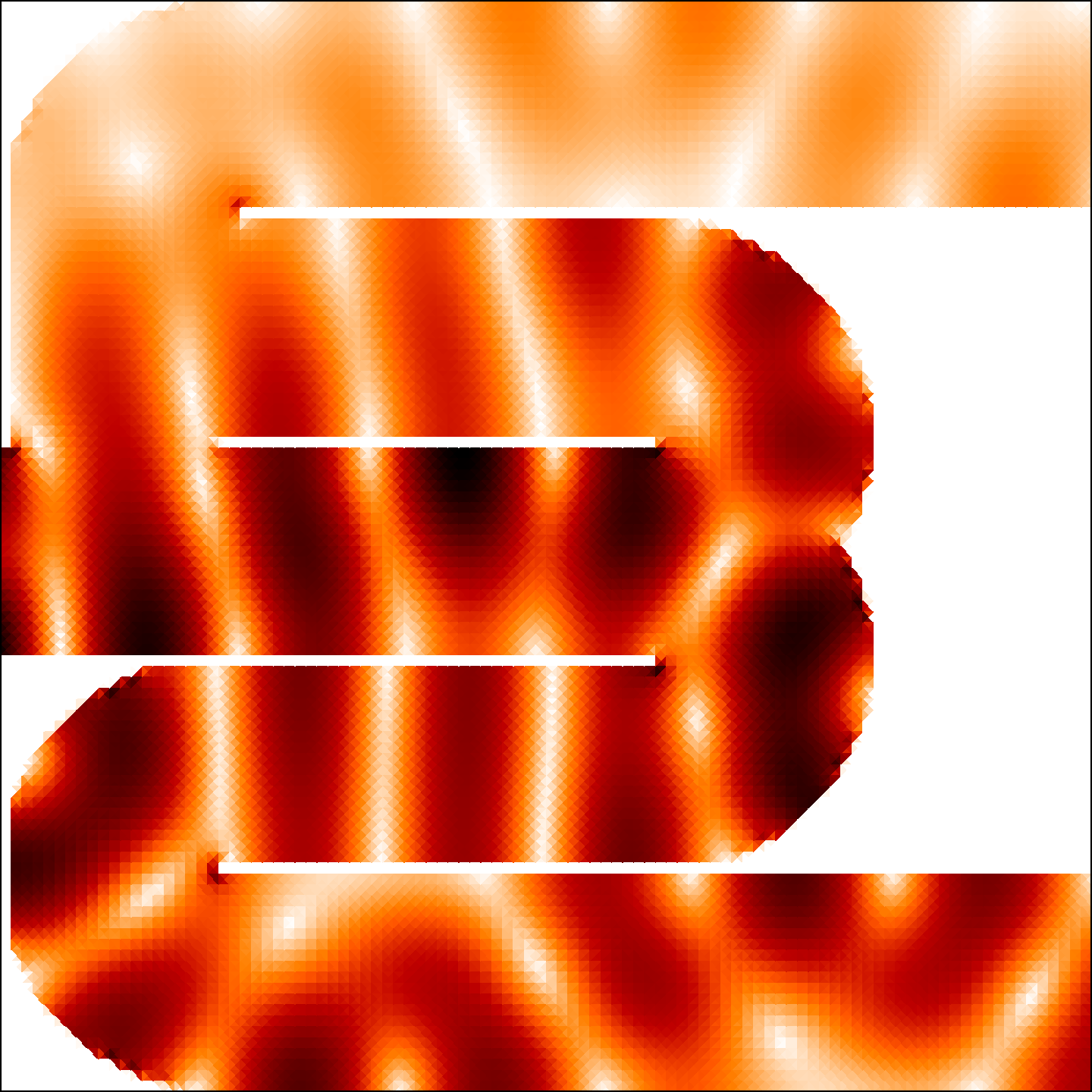}
\end{tabular}
\end{center}
\caption{Example solutions for f=186\,MHz, f=561\,MHz and f=1\,GHz for the first and second geometry. Plotted is $|\mathrm{Re}(E)|$. Script: maxwell\_create\_solutions.py}
\label{buhr_andreas:fig:solutions}
\end{figure}
To get an impression of the solutions, some example solutions are plotted in Fig.~\ref{buhr_andreas:fig:solutions}.
\subsection{Global Properties of Example}
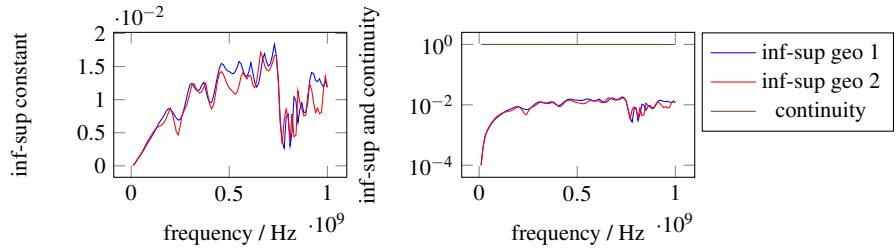
\begin{figure}
\begin{tikzpicture}
\begin{axis}[
    xlabel=frequency / Hz,
    ylabel=inf-sup constant,
    width=0.40\textwidth,
    height=0.3\textwidth,
    no markers,
  ]
  \addplot table{infsups0.txt};
  \addplot table{infsups1.txt};
\end{axis}
\end{tikzpicture}
\begin{tikzpicture}
\begin{semilogyaxis}[
    xlabel=frequency / Hz,
    ylabel=inf-sup and continuity,
    width=0.40\textwidth,
    height=0.3\textwidth,
    no markers,
    legend pos=outer north east,
  ]
  \addplot table{infsups0.txt};
  \addplot table{infsups1.txt};
  \addplot table{continuities0.txt};
  \legend{inf-sup geo 1, inf-sup geo 2, continuity};
\end{semilogyaxis}
\end{tikzpicture}
\caption{inf-sup and continuity constant of bilinear form. Linear and logarithmic. Script: maxwell\_calculate\_infsup.py}
\label{buhr_andreas:fig:infsupcont}
\end{figure}

Before analyzing the behavior of the localized model reduction, we discuss some
properties of the full model. For its stability, its continuity constant $\gamma$ and reduced inf-sup constants $\widetilde \beta$ are the primary concern. They guarantee existence and uniqueness of the solution and their quotient enters the best-approximation inequality
\begin{equation}
\|u - \widetilde u\| \leq \left( 1 + \frac{\gamma}{\widetilde \beta}\right) \inf_{v \in \widetilde V_h} \|u - v\|
\end{equation}
where $u$ is the solution in $V_h$.
Due to the construction of the norm, the continuity constant cannot be larger than one, and numerics indicate that it is usually one (Fig. \ref{buhr_andreas:fig:infsupcont}). The inf-sup constant approaches zero when the frequency goes to zero. This is the well known low frequency instability of this formulation. There are remedies to this problem, but they are not considered here. The order of magnitude of the inf-sup constant is around $10^{-2}$: Due to the Robin boundaries, the problem is stable. With Dirichlet boundaries only, the inf-sup constant would drop to zero at several frequencies.
There are two drops in the inf-sup constant at ca.~770 MHz and 810 MHz. These correspond to resonances in the structure which arise when half a wavelength is the width of a channel ($\lambda / 2 \approx 1/5$).
\begin{figure}
\begin{center}
\begin{tikzpicture}
\begin{semilogyaxis}[
    xlabel=basis size,
    ylabel=max rel. projection error,
    ymax=1e2,
    width=0.45\textwidth,
    height=0.3\textwidth,
    no markers,
    legend pos=outer north east,
  ]
  \addplot table{n_width0.txt};
  \addplot table{n_width1.txt};
  \legend{geometry 1, geometry 2};
\end{semilogyaxis}
\end{tikzpicture}
\end{center}
\caption{
Error when approximating the solution set for all $f \in \varXi$ with an n-dimensional basis obtained by greedy approximation of this set. This is an upper bound for the Kolmogorov n-width.
Script: maxwell\_global\_n\_width.py
}
\label{buhr_andreas:fig:global_n_width}
\end{figure}
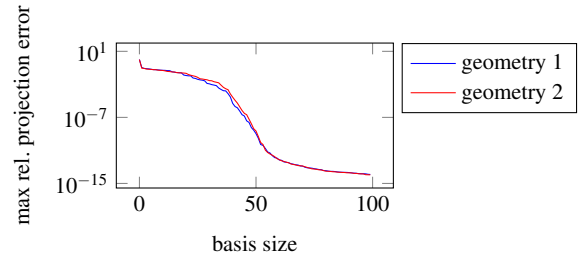

The most important question for the applicability of any reduced basis method is: is the system reducible at all, i.e. can the solution manifold be approximated with a low dimensional solution space? The best possible answer to this question is the Kolmogorov n-width. We measured the approximation error when approximating the solution manifold with a basis generated by a greedy algorithm. The approximation is done by orthogonal projection onto the basis. This error is an upper bound to the Kolmogorov n-width. Already with a basis size of 38, a relative error of $10^{-4}$ can be achieved, see Fig.~\ref{buhr_andreas:fig:global_n_width}. So this problem is well suited for reduced basis methods.

\subsection{Properties of Localized Spaces}
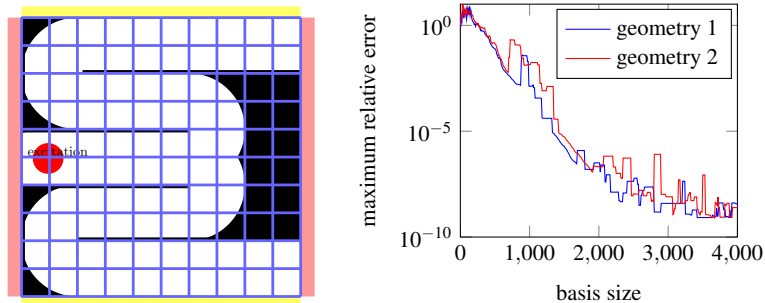
\begin{figure}
\begin{center}
\includegraphics[width=0.35\textwidth]{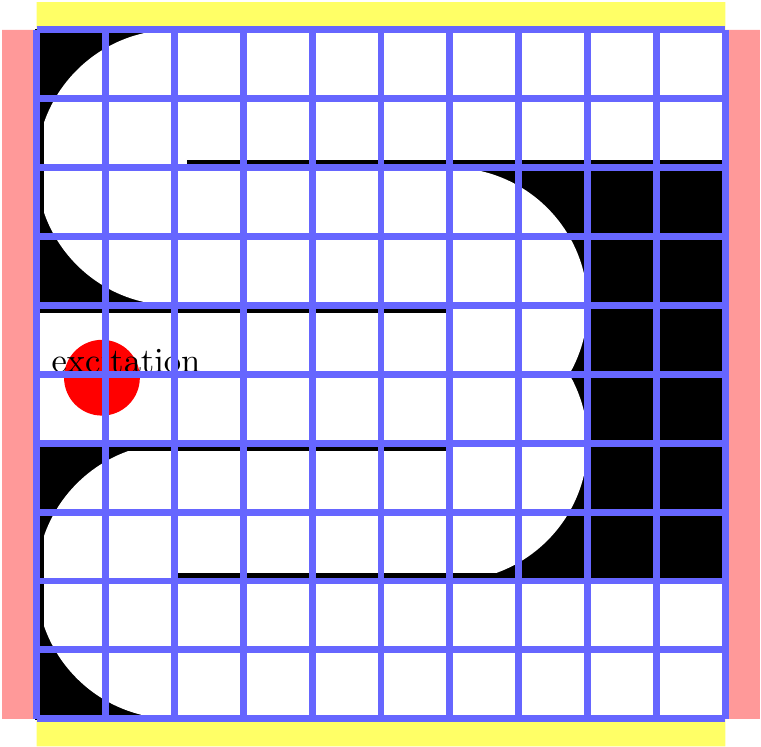}
\hspace{10pt}
\begin{tikzpicture}
\begin{semilogyaxis}[
    xlabel=basis size,
    ylabel=maximum relative error,
    ymin=1e-10,
    ymax=1e1,
    width=0.45\textwidth,
    height=0.4\textwidth,
    xmin=0,
    xmax=4000,
    no markers
  ]
  \addplot table[x index=0, y index=4]{localized0.txt};
  \addplot table[x index=0, y index=4]{localized1.txt};
  \legend{geometry 1, geometry 2};
\end{semilogyaxis}
\end{tikzpicture}
\end{center}
\caption{
Left: Domain decomposition used.
Right: Maximum error when solving with a localized basis, generated by global solves. Script maxwell\_local\_n\_width.py.
}
\label{buhr_andreas:fig:localized_globalsolve}
\end{figure}
The next question is: how much do we lose by localization? Using basis vectors with limited support, one needs a larger total number of basis functions. To quantify this, we compare the errors with global approximation from the previous section with the error obtained when solving with a localized basis, using the best localized basis we can generate. We use a 10 x 10 domain decomposition (see Fig. \ref{buhr_andreas:fig:localized_globalsolve} left) and the space decomposition introduced in Section \ref{buhr_andreas:sec:space_decomposition}. To construct the best possible basis, we solve the full problem for all parameters in the training set. For each local subspace, we apply the corresponding projection operator $P_{\{i\}}$ / $P_{\{i,j\}}$ to all global solutions and subsequently generate a basis for these local parts of global solutions using a greedy procedure. The error when solving in the resulting reduced space is depicted in Fig.~\ref{buhr_andreas:fig:localized_globalsolve}, right. Much more basis vectors are needed, compared to the global reduced basis approach. However the reduction in comparison to the full model (60200 dofs) is still significant and in contrast to standard reduced basis methods, the reduced system matrix is not dense but block-sparse. For a relative error of $10^{-4}$, 1080 basis vectors are necessary.

\begin{figure}
\begin{tikzpicture}
\begin{semilogyaxis}[
    xlabel=basis size,
    ylabel=maximum relative error,
    width=0.42\textwidth,
    height=0.38\textwidth,
    xmin=0,
    xmax=2500
  ]
  \addplot+[
    no markers
  ]table[x index=0, y index=4]{localized0.txt};
\end{semilogyaxis}
\end{tikzpicture}
\begin{tikzpicture}
\begin{semilogyaxis}[
    xlabel=basis size,
    ylabel=reduced inf-sup constant,
    width=0.42\textwidth,
    height=0.38\textwidth,
    xmin=0,
    xmax=2500,
    no markers,
    legend pos=outer north east,
  ]
  \addplot table[
    x expr=\coordindex * 10 + 1,
    y index=16,
  ]{infsup_constant_0.txt};
  \addplot table[
    x expr=\coordindex * 10 + 1,
    y index=33,
  ]{infsup_constant_0.txt};
  \addplot table[
    x expr=\coordindex * 10 + 1,
    y index=69,
  ]{infsup_constant_0.txt};
  \addplot table[
    x expr=\coordindex * 10 + 1,
    y index=79,
  ]{infsup_constant_0.txt};
  \addplot table[
    x expr=\coordindex * 10 + 1,
    y index=83,
  ]{infsup_constant_0.txt};
  \addplot table[
    x expr=\coordindex * 10 + 1,
    y index=89,
  ]{infsup_constant_0.txt};
  \legend{170 MHz, 340 MHz, 700 MHz, 800 MHz, 840 MHz, 900 MHz}
\end{semilogyaxis}
\end{tikzpicture}
\caption{Comparison of maximum error over all frequencies with inf-sup constant of reduced system at selected frequencies for geometry 1. Basis generated by global solves. Increased error and reduced inf-sup constant around basis size of 900. Script: maxwell\_infsup\_during\_reduction.py}
\label{buhr_andreas:fig:infsup_drop}
\end{figure}
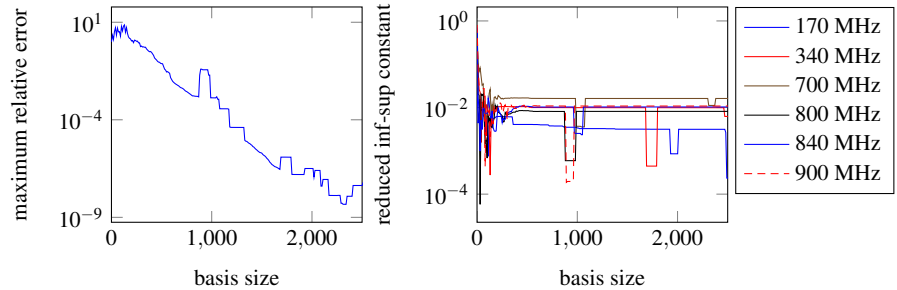

In Fig.~\ref{buhr_andreas:fig:localized_globalsolve} the error is observed to jump occasionally. This is due to the instability of a Galerkin projection of an inf-sup stable problem. While coercivity is retained during Galerkin projection, inf-sup stability is not. While the inf-sup constant of the reduced system is observed to be the same as the inf-sup constant of the full system most of the time, sometimes it drops. This is depicted in Fig.~\ref{buhr_andreas:fig:infsup_drop}.
For a stable reduction, a different test space is necessary. 
However, the application of the known approaches such as
\cite{carlberg2011efficient,dahmen2014double} 
to the localized setting is not straightforward. The development of stable test spaces in the localized setting is beyond the scope of this publication.
\subsection{Properties of Training}
\begin{figure}
\begin{tikzpicture}
\begin{semilogyaxis}[
    xlabel=basis size,
    ylabel=error,
    width=0.48\textwidth,
    height=0.45\textwidth,
    xmin=0,
    xmax=3000,
    ymin=1e-10,
    ymax=1e3,
    no markers,
    title=geometry 1,
  ]
  \addplot table[x index=0, y index=4]{localized0.txt};
  \addplot table[x index=0, y index=4]{maxwell_training_benchmark0.txt};
  \legend{global solves, local training}
\end{semilogyaxis}
\end{tikzpicture}
\begin{tikzpicture}
\begin{semilogyaxis}[
    xlabel=basis size,
    ylabel=error,
    width=0.48\textwidth,
    height=0.45\textwidth,
    xmin=0,
    xmax=3000,
    ymin=1e-10,
    ymax=1e3,
    no markers,
    title=geometry 2,
  ]
  \addplot table[x index=0, y index=4]{localized1.txt};
  \addplot table[x index=0, y index=4]{maxwell_training_benchmark1.txt};
  \legend{global solves, local training}
\end{semilogyaxis}
\end{tikzpicture}
\caption{Maximum error over all frequencies for both geometries. Basis generated by global solves vs.\ basis generated by local training. Script: maxwell\_training\_benchmark.py}
\label{buhr_andreas:fig:training_error}
\end{figure}
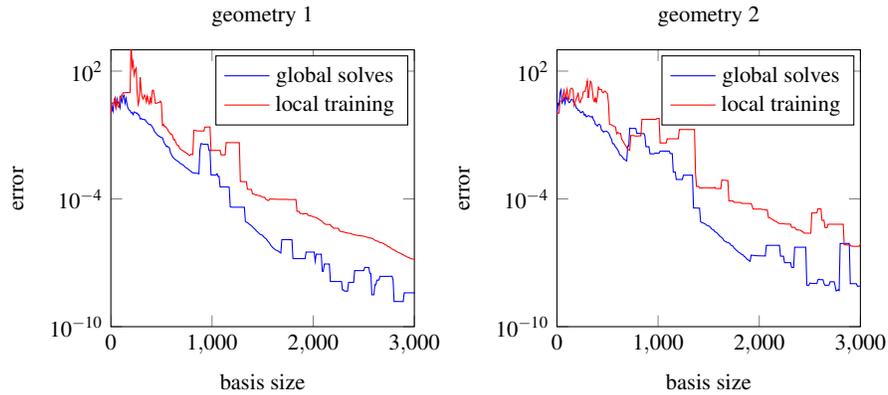

Local basis vectors should be generated using the localized training described in Section \ref{buhr_andreas:sec:training} and in \cite{BEOR15}.
To judge on the quality of these basis vectors, we compare the error obtained using these basis vectors with the error obtained with local basis vectors generated by global solves. The local basis vectors generated by global solves are the reference: These are the best localized basis we can generate. The results for both geometries are depicted in Fig.~\ref{buhr_andreas:fig:training_error}.
While the error decreases more slowly, we still have reasonable basis sizes with training. For a relative error of $10^{-4}$, 1280 basis vectors are necessary for geometry 1 and 1380 are necessary for geometry 2.
\subsection{Application to Local Geometry Change}
\begin{figure}
\begin{center}
\includegraphics[width=0.75\textwidth]{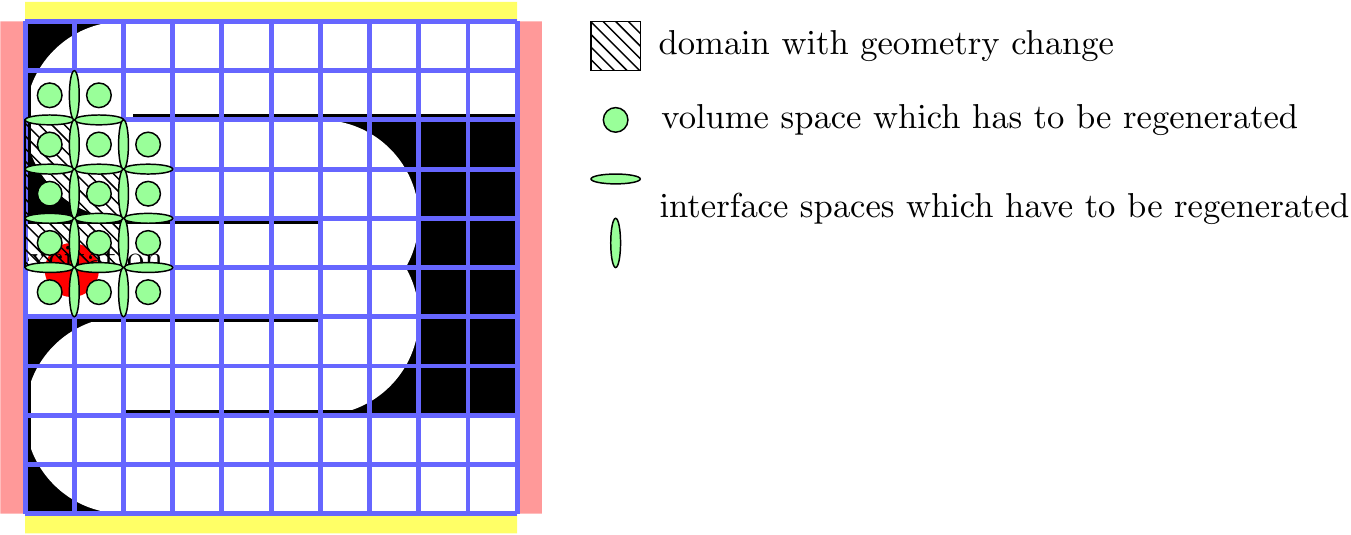}
\end{center}
\caption{Impact of geometry change: 5 domains contain changes, 14 domain spaces and 20 interface spaces have to be regenerated.}
\label{buhr_andreas:fig:example_change}
\end{figure}

\begin{figure}
\begin{center}
geo 1:
\begin{minipage}[c]{0.3\textwidth}
\begin{center}
\includegraphics[width=\textwidth]{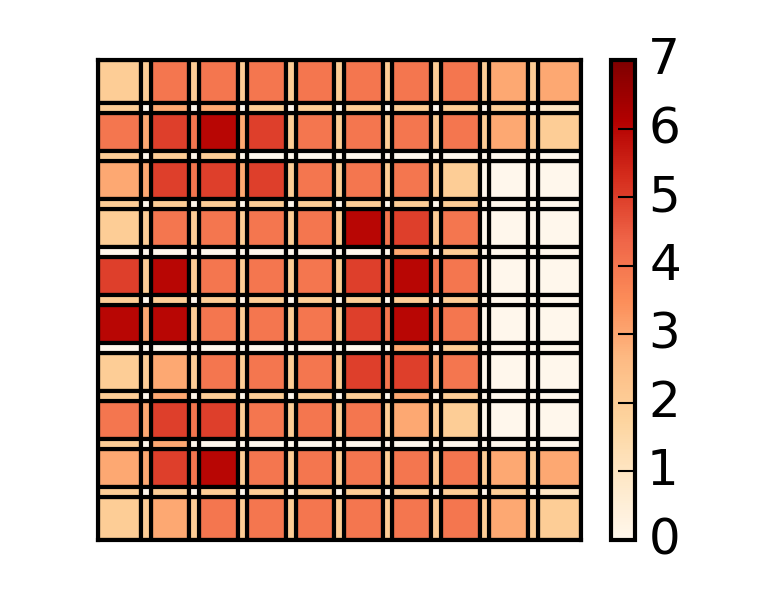}
\end{center}
\end{minipage}
\hspace{20pt}
geo 2:
\begin{minipage}[c]{0.3\textwidth}
\begin{center}
\includegraphics[width=\textwidth]{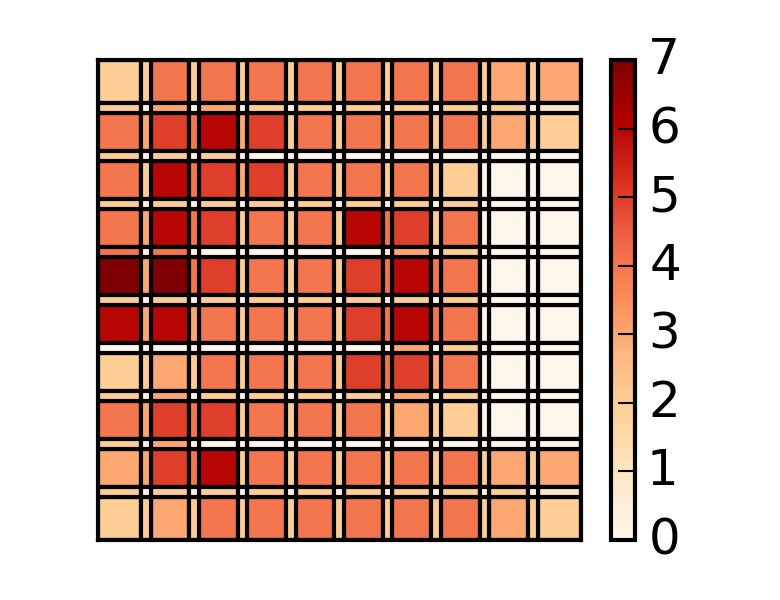}
\end{center}
\end{minipage}
\end{center}
\caption{Basis size distribution. Script: postprocessing\_draw\_basis\_sizes\_maxwell\_geochange.py}
\label{buhr_andreas:fig:basis_sizes}
\end{figure}

If we work with a relative error of 5\%, a basis of size 650 is sufficient for the first geometry and size 675 for the second.
After the geometry change, the local reduced spaces which have no change in their training domain can be reused.
Instead of solving the full system with 60200 degrees of freedom, 
the following effort is necessary per frequency point (see also Fig.~\ref{buhr_andreas:fig:example_change}).
Because the runtime is dominated by matrix factorizations, we focus on these. 
\begin{itemize}
\item
14 factorizations of local problems with 5340 dofs (volume training)
\item
20 factorizations of local problems with 3550 dofs (interface training)
\item
1 factorization of global reduced problem with 675 dofs (global solve)
\end{itemize}
The error between the reduced solution and the full solution in this case is 4.3\%. Script to reproduce: experiment\_maxwell\_geochange.py.
The spacial distribution of basis sizes is depicted in Fig.~\ref{buhr_andreas:fig:basis_sizes}.
\section{Conclusion}
\label{buhr_andreas:sec:conclusion}
ArbiLoMod was applied to the non-coercive problem of 2D Maxwell's equations in $H(\mathrm{curl})$. Its localized training generates a basis of good quality. A reduced model with little error for the full problem can be generated using only local solves, which can easily be parallelized. After localized changes to the model, only in the changed region the localized bases have to be regenerated. All other bases can be reused, which results in large 
computational savings compared to a simulation from scratch. 
The amount of savings is very dependent of the model and the changes which are made. 
A thorough analysis of the computational savings is subject to future work, as is the adaptation of ArbiLoMod's localized a-posteriori error estimator and online enrichment to this problem as well as the instability of the Galerkin projection.
\subsubsection*{Acknowledgments}
Andreas Buhr was supported by CST - Computer Simulation Technology AG.
Stephan Rave was supported by the German Federal Ministry of Education and Research (BMBF) under contract number 05M13PMA.
%bib
\bibliographystyle{plain}
\bibliography{buhr_andreas.bib}
\end{document}